\documentclass{article}
\usepackage{amsfonts}
\usepackage{amsmath}
\usepackage{amsthm}
\usepackage[english]{babel}
\usepackage{graphicx}
\input epsf

\newtheorem{theorem}{Theorem}
\newtheorem{lemma}{Lemma}
\newtheorem{cor}{Corollary}
\newtheorem{note}{Remark}

\begin{document}
\title{In-Degree and PageRank of Web pages:\\
Why do they follow similar power laws?}

\author{N. Litvak\thanks{The work is supported by NWO Meervoud grant no.~632.002.401}, W.R.W. Scheinhardt and Y. Volkovich\thanks{corresponding
author}}

\date{\small{University of Twente, Dept. of Applied Mathematics,\\
P.O.~Box~217, 7500AE~Enschede, The Netherlands; \endgraf e-mail:
\{n.litvak, w.r.w.scheinhardt, y.volkovich\}@ewi.utwente.nl}}

\maketitle

\begin{abstract}
The PageRank is a popularity measure designed by Google to rank
Web pages. Experiments confirm that the PageRank obeys a `power
law' with the same exponent as the In-Degree. This paper presents
a novel mathematical model that explains this phenomenon. The
relation between the PageRank and In-Degree is modelled through a
stochastic equation, which is  inspired by the original definition
of the PageRank, and is analogous to the well-known distributional
identity  for the busy period in the $M/G/1$ queue. Further, we
employ the theory of regular variation and Tauberian theorems to
analytically prove that the tail behavior of the PageRank and the
In-Degree differ only by a multiplicative factor, for which we
derive a closed-form expression. Our analytical results are in
good agreement with experimental data.
\end{abstract}

\section{Introduction}

The notion of {\em PageRank} was introduced by Google in order to
numerically characterize popularity of Web pages. The original
description of the PageRank presented in \cite{BrinPage} is as
follows:
\begin{equation}
PR(i) = c \sum_{j \to i} \frac{1}{d_{j}} PR(j) + (1-c), \label{e1}
\end{equation}
where $PR(i)$ is the PageRank of page $i$, $d_j$ is the number of
outgoing links of page $j$, the sum is taken over all pages $j$
that link to page $i$, and $c$ is the ``damping factor", which is
some constant between 0 and 1. From this equation it is clear that
the PageRank of a page depends on the number of pages that link to
it and the importance (i.e. PageRanks) of these pages.

In this paper we study the relation between the probability
distribution of the PageRank and the {\em In-Degree} of a randomly
selected Web page, where the In-Degree denotes simply the number
of incoming hyperlinks of a Web page. Pandurangan et
al.~\cite{Pandurangan} observed that the probability distributions
of the PageRank and the In-Degree for Web data have a similar
asymptotic behavior, or, more precisely, they seem to follow power
laws with the same exponent. Loosely speaking, a `power law' with
exponent $\alpha$ means that the probability that the random
variable takes some large value $x$ is proportional to
$x^{-\alpha}$. For the PageRank and the In-Degree distribution,
the exponent $\alpha$ is approximately $2.1$.

Recent extensive experiments by Donato et al.~\cite{Donato} and
Fortunato et al.~\cite{Fortunato} confirmed the similarity in tail
behavior observed in~\cite{Pandurangan}. Becchetti and
\\Castillo~\cite{Becchetti} extensively investigated the influence
of the damping factor $c$ on the power law behavior of the
PageRank. They have shown that the PageRank of the top $10\%$ of
the nodes always follows a power law with the same exponent
independent of the value of the damping factor. Our own
experiments based on Web data from \cite{Stanford} are also in
agreement with~\cite{Pandurangan} (see Figure~1).

Obviously, equation~(\ref{e1}) suggests that the PageRank and the
In-Degree are intimately related, but this formula by itself does
not explain the observed similarity in tail behavior. Furthermore,
the linear algebra methods that have been commonly used in the
PageRank literature~\cite{Berkhin,Langville} and proved very
successful for designing efficient computational methods, seem to
be insufficient for modelling and analyzing the  asymptotic
properties of the PageRank distribution.

The goal of our paper is to provide mathematical evidence for the
power-law behavior of the PageRank and its relation to the
In-Degree distribution. We propose a stochastic model that aims to
explain this phenomenon. Our approach is inspired by the
techniques from applied probability and stochastic operations
research. The relation between the PageRank and the In-Degree is
modelled through a distributional identity which is analogous to
the equation for the busy period in the M/G/1 queue (see
e.g.~\cite{Robert}). Further, we analyze our model using the
approach employed in \cite{MeyerTeugels} for studying the tail
behavior of the busy period in case the service times are
regularly varying random variables. This fits in our research
because regular variation is in fact a generalization of the power
law, and it has been widely used in queueing theory to model
self-similarity, long-range dependence and heavy
tails~\cite{Zwart}. Thus, we use the notion of regular variation
to model the power law distribution of the In-Degree. For the sake
of completeness, in Section~\ref{section_prelim}, we will
introduce regularly varying random variables and describe their
basic properties.

To obtain the tail behavior of the PageRank in our model, we use
Laplace-Stieltjes transforms and apply Tauberian theorems
presented in the well-known paper by Bingham and Doney
\cite{BinghamDoney}, see also Theorem~8.1.6 in
\cite{BinghamGoldie}. Moreover, our analysis allows to explicitly
derive the constant multiplicative factor that quantifies the
difference between the PageRank and the In-Degree tail behavior.
Our analytical results show a remarkable agreement with real Web
data.

We believe that our approach is extremely promising for analyzing
the PageRank distribution and solving other problems related to
the structural properties of the Web. At the end of this paper, we
will briefly mention other possibilities for probabilistic
analysis of the PageRank distribution. In particular, we provide
experimental results for Growing Networks~\cite{Barabasi}, and
draw a parallel between the recent studies \cite{Avrachenkov2006,
Fortunato2} on the PageRank behavior in this class of graph models
and our present work.

\section{Preliminaries}
\label{section_prelim}

This section describes important properties of regularly varying
random variables. We follow definitions and notations by Bingham
and Doney~\cite{BinghamDoney}, Meyer and
Teugels~\cite{MeyerTeugels} and Zwart~\cite{Zwart}. More
comprehensive details can be found in \cite{BinghamGoldie}.

We say that a function $V(x)$ is {\em regularly
varying} of index $\alpha \in \mathbb{R}$ if for every $t>0$,
\begin{equation*}
\frac{V(tx)}{V(x)}\to t^{\alpha}\quad\mbox{ as }\quad x\to \infty.
\end{equation*}
If $\alpha=0$, then $V$ is called {\em slowly varying}. Slowly
varying functions are usually denoted by $L$: for every $t>0$,
\begin{equation*}
\frac{L(tx)}{L(x)}\to 1\quad\mbox{ as }\quad x\to \infty.
\end{equation*}
Then, a function $V(x)$ is regularly varying if and only if it can
be written in the form
\begin{equation*}\label{e01}
    V(x)=x^{\alpha}L(x),
\end{equation*}
for some slowly varying $L(x)$.

The following lemma provides a useful bound for slowly varying
functions.
\begin{lemma}
\label{PotterBound} (Potter bounds) Let $L$ be a slowly varying
function. Then, for any fixed $A>1, \delta>0$ there exists a
finite constant $K>1$ such that for all $x_{1}, x_{2}>K,$
\begin{equation*}
\frac{L(x_{1})}{L(x_{2})}\leq A \max \left\{
\left(\frac{x_{1}}{x_{2}}\right)^{\delta},\left(\frac{x_{1}}{x_{2}}\right)^{-\delta}
\right\}.
\end{equation*}
\end{lemma}

In probability theory a random variable $X$ is said to be {\em
regularly varying} with index (or exponent) $\alpha$ if its
distribution function $F_{X}$ is such that
\begin{equation*}
\bar{F}_{X}(x):=1-F_{X}(x)\sim x^{-\alpha}L(x) \quad\mbox{ as
}\quad x\to \infty,\label{e02}
\end{equation*}
for some positive slowly varying function $L(x)$. Here, as in the
remainder of this paper, the notation $a(x)\sim b(x)$ means that
$a(x)/b(x) \to 1$.

We denote the Laplace-Stieltjes transform of $X$ by $f$ and the
$n$th moment $\int_{0}^{\infty}x^{n}dF(x)$ by the corresponding
letter $\mu_{n}$. The successive moments of $F$ can be obtained by
expanding $f$ in a series at $s=0$. More precisely, we have the
following.

\begin{lemma}
\label{AppLem}The $n$th moment of $X$ is finite if and only if
there exist numbers $\mu_{0}=1$ and $\mu_{1},...,\mu_{n},$ such that
$$
f(s)-\sum_{i=0}^{n}\frac{\mu_{i}}{i!}(-s)^{i}=o(s^{n})\mbox{ as } s\to 0.
$$
\end{lemma}

If $\mu_{n}<\infty$ then we introduce the notation $(n\in
\mathbb{N})$
\begin{equation}
f_{n}(s)=(-1)^{n+1}\left(f(s)-\sum_{i=0}^{n}\frac{\mu_{i}}{i!}(-s^{i})\right).
\label{ndef}
\end{equation}

\begin{note}\label{n1}
It follows from Lemma \ref{AppLem} that the $n$th moment of $X$ is
finite if and only if there exist numbers
$\mu_{0}=1$ and $\mu_{1},...,\mu_{n}$ such that $f_{n}(s)=o(s^{n})\mbox{
as }s\to 0.$
\end{note}

The following theorem establishes the relation between asymptotic behavior of regularly varying distribution and its Laplace-Stieltjes transform. This result will play an essential role in our analysis.

\begin{theorem}
\label{AppTh} (Tauberian Theorem) If $n\in \mathbb{N}$,
$\mu_n<\infty$, $\alpha=n+\beta$, $\beta \in (0,1)$, then the
following are equivalent \item[(i)]$f_{n}(s)\sim
(-1)^{\alpha}\Gamma(1-\alpha)
s^{\alpha}L(\frac{1}{s})\quad\mbox{as}\quad s\to 0,$ \item[(ii)]
$1-F(x)\sim x^{-\alpha}L(x)\quad\mbox{as}\quad x\to\infty$.
\end{theorem}

Here and in the remainder of the paper we use the letter $\alpha$
to denote the index of a complementary distribution function
rather than a density. The power law exponent of the In-Degree in
the Web graph then becomes $1.1$ rather than $2.1$.

\section{The model}
\label{model}

In this section we introduce a model that describes the relation
between the PageRank and the In-Degree distributions in the form of
a stochastic equation. This model naturally follows from the
definition of the PageRank (\ref{e1}), and is analytically
tractable for obtaining the asymptotic behavior of the PageRank.

\subsection{Relation between In-Degree and PageRank}
\label{section_pr}

Our goal now is to describe the relation between the PageRank and
the In-Degree. To this end,
we keep equation (\ref{e1}) almost unchanged but we make several
assumptions. First, let $R$ be the PageRank of a randomly chosen
page. We treat $R$ simply as a random variable whose distribution
we want to determine. Second, we assume that the number of
outgoing links $d$ is the same for each page. Then $R$ satisfies a
distributional identity
\begin{equation}
R\stackrel{d}=c\sum_{j=1}^{M}\frac{1}{d}R_{j}+(1-c),  \label{e2a}
\end{equation}
where $M$ is the In-Degree of the considered random page.

We now make the assumption that the $R_j$'s are independent and
have the same distribution as $R$ itself. We note that the independence
assumption is obviously not true in general.
However, it is also not the case that the PageRank values of the pages linking to the same page~$i$ are directly related, so we may assume independence in this
study.

The novelty of our approach is that we treat the PageRank as a
random variable which solves a certain stochastic equation.
However, this approach is quite natural if our goal is to explain
the `power law' behavior of the PageRank because the `power law'
is merely a description of a certain class of probability
distributions. In fact, this point of view is in line with
empirical results by Pandurangan et al.~\cite{Pandurangan} and
other authors who consistently present the (log-log) {\it
histogram} of the PageRank.

One of  the nice features of stochastic equation (\ref{e2a}) is
that it has the same form as the original formula (\ref{e1}).
Thus, we may hope that our model correctly describes the relation
between the In-Degree and the PageRank. This is easy to verify in
the extreme (unrealistic) case when all pages have the same
In-Degree $d$. In this situation, the PageRanks of all pages are
equal, and it is easy to verify that $R\equiv 1$ constitutes the
unique solution of (\ref{e2a}).

\subsection{In-Degree Distribution}
\label{section_nt}

It is well-known that the In-Degree of Web pages follows a power
law. For our analysis however we need a
more formal description of this random variable, thus, we suggest
to employ the theory of regular variation. We model the In-Degree
of a randomly chosen page as a nonnegative, integer, regularly
varying random variable, which is distributed as $N(T)$, where $T$
is regularly varying with index $\alpha$ and $N(t)$ is the number
of Poisson arrivals on the time interval $[0,t]$. Without loss of
generality, we assume that the rate of the Poisson process is
equal to $1$.

The advantage of this construction is that we do not need to
impose any restrictions on $T$  and at the same time ensure that
the In-Degree is integer. We claim that the random variable $N(T)$
will also be regularly varying with the same
 index as $T$, or, more informally, $N(T)$ follows a power
law with the same exponent. Thus, we can think of $N(T)$ as the
In-Degree of a random Web page. For the sake of completeness we
present the formal statement and its proof in the remainder of
this section.

Let $F_{T}$ and $F_{N(T)}$, $f$ and $\phi$ be the distribution
functions and the Laplace-Stieltjes transforms of $T$ and $N(T)$,
respectively. Since the random variable $T$ is regularly varying,
we have by definition
\begin{equation}
1-F_{T}(x)\sim x^{-\alpha}L(x)\;\mbox{ as } x\to\infty,
\label{ein1}
\end{equation}
where $L(x)$ is some slowly varying function. Then we will claim that for
$N(T)$ the following also holds:
\begin{equation}
1-F_{N(T)}(x)\sim x^{-\alpha}L(x)\;\mbox{ as }
x\to\infty.\label{ein2}
\end{equation}

To prove this statement we use the Tauberian theorem (Theorem
\ref{AppTh}). In order to satisfy the conditions of this theorem,
we should first verify whether the corresponding moments of $T$
and $N(T)$ always exist together. Assuming that $\mathbb{E}T=d$ we
immediately get $\mathbb{E}N(T)= d$. Next, consider the generating
function of $N(T)$,
\begin{equation}
\mathbb{G}_{N(T)}(s):=\mathbb{E}s^{N(T)}=\int_{0}^{\infty}\mathbb{E}s^{N(t)}dF_{T}(t)=\int_{0}^{\infty}e^{-t(1-s)}dF_{T}(t)=f(1-s),
\label{e14}
\end{equation}
from which we derive the Laplace-Stieltjes transform of $N(T)$ in
terms of the Laplace-Stieltjes transform of $T$:
\begin{equation*}
\phi(w)=\mathbb{E}e^{-wN(T)}=f(1-e^{-w}).
\end{equation*}

Now, denote by $\mu_{1},\ldots,\mu_{n}$ and $\xi_{1},\ldots,\xi_{n}$
the first $n$ moments of $T$ and $N(T)$, respectively, and define $\mu_0=\xi_0=1$. Then we can formulate the
next lemma.
\begin{lemma}
The following are equivalent
\begin{itemize}
\item[(i)]  $\mu_{n}<\infty,$
\item[(ii)] $\xi_{n}<\infty.$ \label{l3}
\end{itemize}
\end{lemma}
\begin{proof}
\item[$(i)\to(ii)$] By Lemma \ref{AppLem} we know that $\mu_{n}<\infty$
if and only if $f(t)$ can be written as
$$
f(t)=\sum_{i=0}^{n}\frac{\mu_{i}}{i!}(-t)^{i}+o(t^{n})\mbox{ as }
t\to 0.
$$
Denote $t(s):=1-e^{-s}$, then $t(s)\to 0$ as $s\to 0$, and we can
substitute
\begin{eqnarray*}
\phi(s)&=&f(1-e^{-s})=\sum_{i=0}^{n}\frac{\mu_{i}}{i!}(-(1-e^{-s}))^{i}+o((1-e^{-s})^{n}) \\
&=&\sum_{i=0}^{n}\frac{\mu_{i}}{i!}(-1)^{i}\left(\sum_{k=1}^{\infty}(-1)^{k+1}\frac{s^k}{k!}\right)^{i}+o(s^{n}),
\end{eqnarray*}
which can be written as
$$
\phi(s)=\sum_{i=0}^{n}\frac{\xi_{i}}{i!}(-s)^{i}+o(s^{n})
$$
for some finite constants $\xi_{0}=1$ and
$\xi_{1},\ldots,\xi_{n}$, that can be expressed in terms of
$\mu_{0}=1$ and $\mu_{1},\ldots,\mu_{n}$. Thus, by
uniqueness of the power series expansion and by Lemma \ref{AppLem}
we have $\xi_{n}<\infty.$

\item[$(ii)\to(i)$] Similarly, $s(t):=-\ln(1-t)\to 0$ as $t\to 0$,
so we obtain
\begin{eqnarray*}
f(t)&=&\phi(-\ln(1-t))=\sum_{i=0}^{n}\frac{\xi_{i}}{i!}\ln^{i}(1-t)+o(\ln^n(1-t))\\
&=&\sum_{i=0}^{n}\frac{\xi_{i}}{i!}\left(-\sum_{k=1}^{\infty}\frac{t^{k}}{k}\right)^{i}+o\left(\left(-\sum_{k=1}^{\infty}\frac{t^{k}}{k}\right)^{n}\right)\\
&=&\sum_{i=0}^{n}\frac{\mu_{i}}{i!}(-t)^{i}+o(t^{n}),
\end{eqnarray*}
for $\mu_{0}=1$ and some $\mu_{1},\ldots,\mu_{n}$
that can be expressed in terms of $\xi_{0}=1$ and
$\xi_{1},\ldots,\xi_{n}$, which similarly implies $\mu_{n}<\infty$.
\end{proof}

\begin{note}
If we define
\begin{equation*}
    f_{n}(s)=(-1)^{n+1}\left(f(s)-\sum_{i=0}^{n}\frac{\mu_{i}}{i!}(-s)^{i}\right)\mbox{ and }
\end{equation*}
\begin{equation*}
    \phi_{n}(s)=(-1)^{n+1}\left(\phi(s)-\sum_{i=0}^{n}\frac{\xi_{i}}{i!}(-s)^{i}\right)
\end{equation*}
as in (\ref{ndef}), then we can reformulate Lemma \ref{l3}
as follows:
\begin{equation*}
    f_{n}(s)=o(s^{n})\quad\mbox{if and only if}\quad\phi_{n}(s)=o(s^{n}).
\end{equation*}
\end{note}

Now, we can use Theorem \ref{AppTh} to prove that (\ref{ein1})
implies (\ref{ein2}). In fact also the reverse holds, as stated in
the next theorem.
\begin{theorem} The following are equivalent
\item[(i)] $\bar{F}_{T}(x)\sim x^{-\alpha}L(x)\quad\mbox{as}\quad
x\to \infty,$ \item[(ii)] $\bar{F}_{N(T)}(x)\sim
x^{-\alpha}L(x)\quad\mbox{as}\quad x\to \infty.$ \label{t3}
\end{theorem}
\begin{proof}
\item[$(i)\to(ii)$] From Theorem~\ref{AppTh} for $T$ we know that
\begin{eqnarray} \nonumber \bar{F}_{T}(x)\sim x^{-\alpha}L(x),\
x\to\infty \quad
\mbox{implies}\\
f_{n}(t)\sim(-1)^{\alpha}\Gamma(1-\alpha)t^{\alpha}L\left(\frac{1}{t}\right)\quad\mbox{as}\quad
t\to0, \label{e3}
\end{eqnarray}
where $\alpha>1$ is not integer and $n$ is the largest integer smaller than $\alpha$.

Since $\phi(s)=f(t)$, by Lemma~\ref{l3} we have
$f_{n}(t)\sim\phi_{n}(s)$, where $t(s)=(1-e^{-s})\sim s$, as $s\to
0$. So, we can obtain from (\ref{e3}) by using Lemma
\ref{PotterBound} that
\begin{equation*}
\phi_{n}(s)\sim(-1)^{\alpha}\Gamma(1-\alpha)s^{\alpha}
L\left(\frac{1}{s}\right).
\end{equation*}
Now we again apply Theorem~\ref{AppTh} to
conclude that
$$
\bar{F}_{N(T)}\sim x^{-\alpha}L(x)\;\mbox{ as }x\to\infty.
$$
\item[$(ii)\to(i)$] Similar to the first part of the proof.
\end{proof}

Thus, our model for the number of incoming links properly
describes an In-Degree distribution that follows a power law with finite expectation and a
 non-integer exponent.

\subsection{The main stochastic equation}

Combining the ideas from Sections~\ref{section_pr} and
\ref{section_nt}, we arrive to the following equation
\begin{equation}
R\stackrel{d}=c\sum_{j=1}^{N(T)}\frac{1}{d}R_{j}+(1-c),
\label{e2}
\end{equation}
where $c\in(0,1)$ is the damping factor, $d\in\{1,2,\ldots\}$ is
the fixed Out-Degree of each page, and $N(T)$ describes the
In-Degree of a randomly chosen page as the number of Poisson
arrivals on a regularly varying time interval $T$. As we
discussed above, stochastic equation (\ref{e2}) adequately
captures several important aspects of the PageRank distribution
and its relation to the In-Degree. Moreover, our model is
completely formalized, and thus we can apply analytical methods in
order to derive the tail behavior of the random variable $R$
representing the PageRank.

Linear stochastic equations like (\ref{e2}) have a long history.
In particular, (\ref{e2}) is similar to the famous equation that
arises in  the theory of branching processes and describes many
real-life phenomena, for instance, the distribution of the busy
period in the $M/G/1$ queue:
\begin{equation*}
B\stackrel{d}{=}\sum_{i=1}^{N(S_1)}B_i+S_1,
\end{equation*}
where $B$ is the distribution of the busy period (the time
interval during which the queue  is non-empty), $S_1$ is the
service time of the customer that initiated the busy period,
$N(S_1)$ is the number of Poisson arrivals during this service
time and the $B_i$'s are independent and distributed as $B$. We
refer to \cite{Robert} and other books on queueing theory for more
details. Also, see Zwart~\cite{Zwart} for an excellent detailed
treatment of  queues with regular variation, and specifically the
busy period problem. We note also that our equation (\ref{e2}) is
a special case in a rich class of stochastic recursive equations
that were discussed in detail in the recent survey by Aldous and
Bandyopadhyay~\cite{Aldous}.

 This concludes the model description. The next step will be to use our model for
providing a rigorous explanation of the indicated connection between the distributions of the In-Degree
and the PageRank.

\section{Analysis}

The idea of our analysis is to write the equation for the
Laplace-Stieltjes Transforms of $T$ and $R$ and then make use of
the Tauberian theorems to prove that $R$ is regularly varying with
the same index as $T$. According to Theorem~\ref{t3}, this will
give us the desired similarity in tail behavior of the PageRank
$R$ and the In-Degree $N(T)$.

As a result of the assumptions from Section~\ref{model}, we can express the Laplace-Stieltjes transform $r(s)$ of the PageRank
distribution $R$ in terms of the probability generating function of $N(T)$ using (\ref{e2}):
\begin{eqnarray*}
r(s)&:=&\mathbb{E}e^{-sR}
=e^{-s(1-c)}\mathbb{E}\exp\left(-s\frac{c}{d}\sum_{i=1}^{N(T)}R_{i}\right)\\
&=&e^{-s(1-c)}\sum_{k=1}^{\infty}\mathbb{E}\exp\left(-s\frac{c%
}{d}\sum_{i=1}^{k}R_{i}\right)\mathbb{P}(N(T)=k)\\
&=&e^{-s(1-c)}\sum_{k=1}^{\infty}\Pi_{i=1}^{k}\mathbb{E}\exp
\left(-s\frac{c}{d}R\right)\mathbb{P}(N(T)=k)\\
&=&e^{-s(1-c)}\sum_{k=1}^{\infty}\left(r\left(s\frac{c}{d}\right)\right)^{k}\mathbb{P}(N(T)=k)\\
&=&e^{-s(1-c)}\mathbb{G}_{N(T)}\left(r\left(s\frac{c}{d}\right)\right).
\end{eqnarray*}
Since, by (\ref{e14}), $\mathbb{G}_{N(T)}(s)=f(1-s),$ we arrive at
\begin{equation}
r(s)=f\left(1-r\left(\frac{c}{d}s\right)\right)e^{-s(1-c)}.
\label{e5}
\end{equation}
It can be shown (e.g. arguing as in
\cite[Section~XIII.4]{Feller2}) that equation (\ref{e5}) has a
unique solution $r(s)$ which is completely monotone and has
$r(0)=1$ if and only if $c/d<1$. This inequality is satisfied for
the typical values $d
> 1$ and $0<c<1$.

As in Section~\ref{section_nt}, we will start the analysis with
providing the correspondence between existence of the $n$-th
moments of $T$ and $R$. We remind that $\mu_{1},\ldots,\mu_{n}$ denote the first $n$ moments of $T$.
 Further, denote the first $n$ moments of $R$ by
$\eta_{1},\ldots,\eta_{n}$, and define
\begin{equation*}
r_{n}(s)=(-1)^{n+1}\left(r(s)-\sum_{k=0}^{n}\frac{\eta_{k}}{k!}(-s^k)\right),
\end{equation*}
as in (\ref{ndef}). Note that taking expectations on both sides of
(\ref{e2}) we easily obtain $\mathbb{E}R=\eta_1=1$. This follows
from the independence of $N(T)$ and the $R_j$'s and the fact that
$\mathbb{E}N(T)=\mathbb{E}T=\mu_1=d$.

The next lemma holds.
\begin{lemma}
\label{l1} The following are equivalent
\begin{itemize}
\item[(i)]
$\mu_{n}<\infty$, \item[(ii)] $\eta_{n}<\infty$.
\end{itemize}
\end{lemma}

\begin{proof}
\item[$(i)\to(ii)$] We use induction, starting from $n=1$ for which both $(i)$ and $(ii)$
are valid. Assume that for $k=1,2,\ldots,n-1$ it has been shown
that $(i)\to(ii)$. We introduce the following notation, to be used throughout this section.
Denote
\begin{eqnarray*}
g(s)&:=&e^{-s(1-c)},\qquad \mbox{and}\\
t(s)&:=&1-r\left(\frac{c}{d}s\right).
\end{eqnarray*}
Then we can write (\ref{e5}) as
\begin{equation}
r(s)=f(t)g(s).
\label{e6}
\end{equation}
We know from $(i)$ that
\begin{align*}
f(t)&=1-dt+\sum_{k=2}^{n}\frac{\mu_{k}(-t)^{k}}{k!}+o(t^n)\\
&=1-d\left(1-r\left(\frac{c}{d}s\right)\right)+\sum_{k=2}^{n}\frac{\mu_{k}(-t)^{k}}{k!}+o(t^{n}).
\end{align*}
Thus, from (\ref{e6}) we obtain
\begin{equation}
r(s)-dg(s)r\left(\frac{c}{d}s\right)=\left(1-d+\sum_{k=2}^{n}\frac{\mu_{k}(-t)^{k}}{k!}+o(t^{n})\right)g(s).\label{e7}
\end{equation}
However, it follows from the induction hypothesis for $n-1$ that
$$
r(s)=1-s+\sum_{k=2}^{n-1}\frac{\eta_{k}}{k!}(-s^k)+o(s^{n-1}),
$$
so we can present $t(s)$ as a sum
$$
t(s)=-\sum_{k=1}^{n-1}\frac{\eta_{k}}{k!}\left(\frac{c}{d}\right)^{k}(-s)^{k}+o(s^{n-1}).
$$
Using this, we can actually find $t^k(s)$:
\begin{equation}
t^{k}(s)=\sum_{i=k}^{n+k-2}\beta_{k,i}s^{i}+o(s^{n+k-2}),\label{e8}
\end{equation}
for $k\ge 1$ and appropriate constants $\beta_{k,i},
i=k,\ldots,k+n-2.$ Thus, we obtain by (\ref{e7}) and (\ref{e8}):
\begin{equation*}
r(s)-d g(s)r\left(\frac{c}{d}s\right)
=\left(\sum_{i=0}^{n}\gamma_{i}(-s)^{i}+o(s^{n})\right)g(s)
\label{e9}
\end{equation*}
for appropriate constants $\gamma_{0},\ldots,\gamma_{n}$. Using
the expansion of $g(s)$, it is not difficult to show that for
appropriate constants $\rho_{0},\ldots,\rho_{n},$ we also have
\begin{equation*}
r(s)-dr\left(\frac{c}{d}s\right)=\sum_{i=0}^{n}\rho_{i}s^{i}+o(s^{n}).\label{e10}
\end{equation*}
In other words, because of the uniqueness of the series expansion, we have
\begin{equation}
\left(r(s)-dr\left(\frac{c}{d}s\right)\right)_{n}=r_{n}(s)-dr_{n}\left(\frac{c}{d}s\right)=o(s^{n}).\label{e11}
\end{equation}
We will now show that this implies $(ii)$, to which end we
consider the partial sums
\begin{eqnarray*}
r_{n}^{N}(s)&=&\sum_{k=0}^{N}d^{k}\left(r_{n}\left(\left(\frac{c}{d}\right)^{k}s\right)-dr_{n}\left(\left(\frac{c}{d}\right)^{k+1}s\right)\right)\\
&=&r_{n}(s)-d^{N+1}r_{n}\left(\left(\frac{c}{d}\right)^{N+1}s\right).
\end{eqnarray*}
Taking the limit as $N\to\infty$, we have for the last term that
\begin{align*}
&\lim_{N\to\infty} d^{N+1}r_{n}\left(\left(\frac{c}{d}\right)^{N+1}s\right)\\
&=
\lim_{N\to\infty} \frac{r_{n}\left(\left(\frac{c}{d}\right)^{N+1}s\right)}
                       {\left(\left(\frac{c}{d}\right)^{N+1}s\right)^{n-1}}
\lim_{N\to\infty}
\left(\frac{c}{d}\right)^{(N+1)(n-2)}s^{n-1}c^{N+1}=0,
\end{align*}
where we used the induction hypothesis $r_{n}(s)=o(s^{n-1})$
together with $n\geq~2,\ 0<c<1$ and $d>1$. It follows that we can
express $r_{n}(s)$ as an infinite sum,
\begin{equation}
\label{e16}
r_{n}(s)=\sum_{k=0}^{\infty}d^{k}\left(r_{n}\left(\left(\frac{c}{d}\right)^{k}s\right)-dr_{n}\left(\left(\frac{c}{d}\right)^{k+1}s\right)\right),
\end{equation}
where we can apply (\ref{e11}) to each of the terms. Further,
by definition of $o(s^{n})$, for all $\varepsilon>0$,
there exists a $\delta=\delta(\varepsilon)$ such that
$\left|r_{n}(s)-dr_{n}\left(\frac{c}{d}s\right)\right|<
\varepsilon s^{n}$ whenever $0<s\le\delta$.
Moreover, for this $\varepsilon$ and $\delta$, and $0<s\le\delta$, we also have
\begin{eqnarray}
\nonumber
|r_{n}(s)|&=&\left|\sum_{k=0}^{\infty}d^{k}\left(r_{n}\left(\left(\frac{c}{d}\right)^{k}s\right)-dr_{n}\left(\left(\frac{c}{d}\right)^{k+1}s\right)\right)\right|\\
\nonumber
&\le&
\sum_{k=0}^{\infty}\left|d^{k}\left(r_{n}\left(\left(\frac{c}{d}\right)^{k}s\right)-dr_{n}\left(\left(\frac{c}{d}\right)^{k+1}s\right)\right)\right|\\
\label{delta}
&<&\sum_{k=0}^{\infty}\varepsilon d^k\left(\frac{c}{d}\right)^{kn}s^{n}=\frac{d^{n-1}}{d^{n-1}-c^{n}}\varepsilon s^{n}.
\end{eqnarray}
Here the second inequality holds because $0<\left(\frac{c}{d}\right)^{k}s\le\delta$ for every
$k\ge0$.
Since for every $\varepsilon_{0}>0$ there exists
$\delta_{0}$ such that $$\left|r_{n}(s)-dr_{n}\left(\frac{c}{d}s\right)\right|<
\frac{d^{n-1}-c^{n}}{d^{n-1}}\varepsilon_{0} s^{n}$$ for $0<s\le\delta_0$, then
according to (\ref{delta}),
we have $|r_{n}(s)|< \varepsilon_{0}s^{n}$ whenever
$0<|s|\le\delta_{0}$, by which we have shown that $r_{n}=o(s^{n}).$

\item[$(ii)\to(i)$]

Assume that there exists a nonnegative random variable $R$
satisfying (\ref{e2}). Then, obviously, $R\ge1-c$. Moreover,
(\ref{e2}) also implies that $R$ is stochastically greater than
$(1-c)\left(\frac{c}{d} N(T)+1\right)$. Hence, the existence of
the $n$-th moment of $R$ ensures the existence of the $n$-th
moment of $N(T)$, which in turn by Lemma~\ref{l3} ensures the
existence of the $n$-th moment of $T$.
\end{proof}

\begin{note} Note that the stochastic inequality $R\stackrel{d}{>}(1-c)\left(\frac{c}{d}N(T)+1\right)$ implies that the tail of the PageRank is at least as heavy as the tail of the In-Degree.
\end{note}

\begin{note}
    Similar  as in Remark~\ref{n1}, we can reformulate Lemma
    \ref{l1} as
    \begin{equation*}
    f_{n}(s)=o(s^{n})\quad \mbox{if and only if}\quad r_{n}(s)=o(s^{n}).
\end{equation*}
\end{note}

From the first part of the proof of Lemma~\ref{l1} we also obtain the next corollary.
\begin{cor} The following holds:
\begin{equation*}
r_{n}(s)-dr_{n}(\frac{c}{d}s)=f_n(t)+O(t^{n+1}).
\end{equation*}
\label{co1}
\end{cor}

\begin{proof}
By definitions of $r_n(s)$, $f_{n}(t)$, $t(s)$ and Lemma~\ref{l1}, it follows from (\ref{e6}) that for fixed $n$,
\begin{align*}
(-1)^{n+1}r_{n}(s)+\sum_{k=0}^{n}\frac{\eta_{k}}{k!}(-s^k)=
\left((-1)^{n+1}f_{n}(t)+1-dt+\sum_{k=2}^{n}\frac{\mu_{k}(-t)^{k}}{k!}\right)g(s)\\
=\left((-1)^{n+1}f_{n}(t)+1-d+d\left((-1)^{n+1}r_{n}\left(\frac{c}{d}s\right)
+\sum_{k=0}^{n}\frac{\eta_{k}}{k!}\left(\frac{c}{d}\right)^{k}(-s)^{k}\right)+\right.\\
\left.+\sum_{k=2}^{n}\frac{\mu_{k}(-t)^{k}}{k!}\right)(1+o(1)).\\
\end{align*}

Because $r_{n}(s)=o(s^{n})$ we can extend (\ref{e8}) for $k\ge 1$
and appropriate constants $\beta_{k,i}, i=k,...,k+n-1$:
\begin{equation*}
t^{k}(s)=\sum_{i=k}^{n+k-1}\beta_{k,i}s^{i}+o(s^{n+k-1}),
\end{equation*}
and rewrite the last equation as
\[(-1)^{n+1}r_{n}(s)+\sum_{k=0}^{n}\frac{\eta_{k}}{k!}(-s^k)\]
\[
=(-1)^{n+1}f_{n}(t)-d(-1)^{n+1}r_{n}\left(\frac{c}{d}s\right)+\sum_{k=0}^{n+1}\tau_{k}s^{k}+
o(s^{n+1}),\]
where $\tau_{0},\ldots,\tau_{n+1}$ are corresponding constants. Now due to the
uniqueness of the series expansion, we can reduce the above formula to
$$
r_{n}(s)=f_{n}(t)+dr_{n}\left(\frac{c}{d}s\right)+(-1)^{n+1}\tau_{n+1}s^{n+1}+o(s^{n+1}).
$$
Then we get:
$$
r_{n}(s)-dr_{n}\left(\frac{c}{d}s\right)=f_n(t)+O(t^{n+1}).
$$
\end{proof}

Now we are ready to explain the similarity between the In-Degree
and the PageRank distributions. The next theorem formalizes this main
statement.
\begin{theorem}
The following are equivalent
\begin{itemize}
\item[(i)] $\bar{F}_{N(T)}(x)\sim
x^{-\alpha}L(x)\quad\mbox{as}\quad x\to \infty,$ \item[(ii)]
$\bar{F}_{R}(x)\sim{\displaystyle\frac{c^{\alpha}}{d^{\alpha}-c^{\alpha}d}}x^{-\alpha}L(x)\quad\mbox{as}\quad
x\to \infty.$
\end{itemize}
\label{t2}
\end{theorem}
\begin{proof}
\item[$(i)\to(ii)$] From $(i)$ and Theorem~\ref{t3} it follows that
\begin{equation}
    \bar{F}_{T}(x)\sim x^{-\alpha}L(x)\quad\mbox{as}\quad x\to \infty.
    \label{e17}
\end{equation}
Theorem~\ref{AppTh} also implies that (\ref{e17}) is equivalent to
$f_{n}(t)\sim(-1)^{\alpha}\Gamma(1-\alpha)t^{\alpha}L\left(\frac{1}{t}\right)$,
where $t(s)\sim (c/d)s,$ as $s\to 0$. Then, by Corollary
{\ref{co1}} we obtain
$$
r_{n}(s)-dr_{n}\left(\frac{c}{d}s\right)\sim(-1)^{n}\Gamma(1-\alpha)\left(\frac{c}{d}\right)^{\alpha}s^{\alpha}L\left(\frac{1}{s}\right)\;\mbox{ as } s\to 0.
$$
Then also for every $k\ge 0$, as $s\to 0$, we have
\begin{eqnarray*}
r_{n}\left(\left(\frac{c}{d}\right)^{k}s\right)-dr_{n}\left(\left(\frac{c}{d}\right)^{k+1}s\right)
&\sim&(-1)^{n}\Gamma(1-\alpha)\left(\frac{c}{d}\right)^{\alpha}\left(\frac{c}{d}\right)^{\alpha
k}s^{\alpha}L\left(\frac{1}{\left(\frac{c}{d}\right)^{k}s}\right)\\
&\sim&
(-1)^{n}\Gamma(1-\alpha)\left(\frac{c}{d}\right)^{\alpha}\left(\frac{c}{d}\right)^{\alpha
k}s^{\alpha}L\left(\frac{1}{s}\right),
\end{eqnarray*}
and from the infinite-sum representation (\ref{e16}) for $r_{n}(s)$, we
directly obtain
\begin{equation*}
r_{n}(s)\sim(-1)^{n}\Gamma(1-\alpha)\frac{d^{\alpha}}{d^{\alpha}-c^{\alpha}d}
\left(\frac{c}{d}\right)^{\alpha}s^{\alpha}L\left(\frac{1}{s}\right)\;\mbox{
as }s\to 0.
\end{equation*}
Now we again apply Theorem~\ref{AppTh}, which leads to $(ii)$.
\item[$(ii)\to(i)$] The proof follows easily from $(ii)$ and
Corollary~\ref{co1}.
\end{proof}

Thus, we have shown that the asymptotic behavior of the PageRank and the In-Degree differ only by the multiplicative factor $\frac{c^{\alpha}}{d^{\alpha}-c^{\alpha}d}$ whereas the power law exponent remains the same. In the next section we will experimentally verify this result.

\section{Numerical Results}
\subsection{Power Law Identification}
    The identification and measuring of power law behavior is not always simple. In this
section we provide a brief overview of techniques that we used to
plot and numerically identify power law distributions.

The standard strategy is to plot a histogram of a quantity on
logarithmic scales to obtain a straight line, which is a typical
feature of the power law. However, this  technique is often not
efficient. In~\cite{Newman}, Newman  clearly illustrated that even
for generated random numbers with a known distribution the noise
in the tail region has a strong influence on the estimation of the
power law parameters. He suggests to plot the fraction of
measurements that are not smaller than a given value, i.e. the
complementary cumulative distribution function $\bar{F}(x) =
P(X\geq x)$ rather than the histogram. The advantage is that we
obtain a less noisy plot. Besides, this idea is consistent with
our analysis in the previous section, which was based on
complementary cumulative distribution functions. We note that if
the distribution of $X$ follows a power law with exponent $\alpha$
so that $\bar{F}(x)\sim Cx^{-\alpha}$, $x\to\infty$, where $C$ is
some constant, then the corresponding histogram has an exponent
$\alpha+1$. Thus, the plot of $\bar{F}(x)$ on logarithmic scales
has a smaller slope than the plot of the histogram.

Computing the correct slope from the observed
data is also not trivial. Goldstein et al. in \cite{Golshtein}, and later Newman in
\cite{Newman}, have proposed to use maximum likelihood estimation,
which provides a more robust estimation of the power law exponent
than the standard least-squares fit method. Thus, we compute the
exponent $\alpha$ using the next formula from~\cite{Newman}:
\begin{equation}\label{e18}
\alpha=1+N\left(\sum_{i=1}^{N}\ln\frac{x_{i}}{x_{min}}\right).
\end{equation}
Here the quantities $x_{i}$, $i=1,\ldots,N$, are the measured values
of $X$, and $x_{min}$ usually corresponds to the
smallest value of $X$ for which the power law behavior is assumed to hold.

In the next sections we will present our experiments on real Web
Data and on a graph that represents a well-known mathematical
model of the Web (Growing Networks). In both cases, for each value
$x$, we plot in log-log scale the number of measurements that are
not smaller than $x$, and we use (\ref{e18}) to obtain the
exponents.

\subsection{Web Data}  \label{sec:WD}
To confirm our results on asymptotic similarity between PageRank
and In-Degree distributions we performed experiments on the public
data of the Stanford Web from \cite{Stanford}. We calculated the
PageRanks over a Web graph with $281903$ nodes (pages) and $\sim
2.3$ million edges (links) using the standard power method (see
e.g. \cite{Langville}).

\begin{figure}[hbt]
            \centering {\epsfxsize=4 in \epsfbox{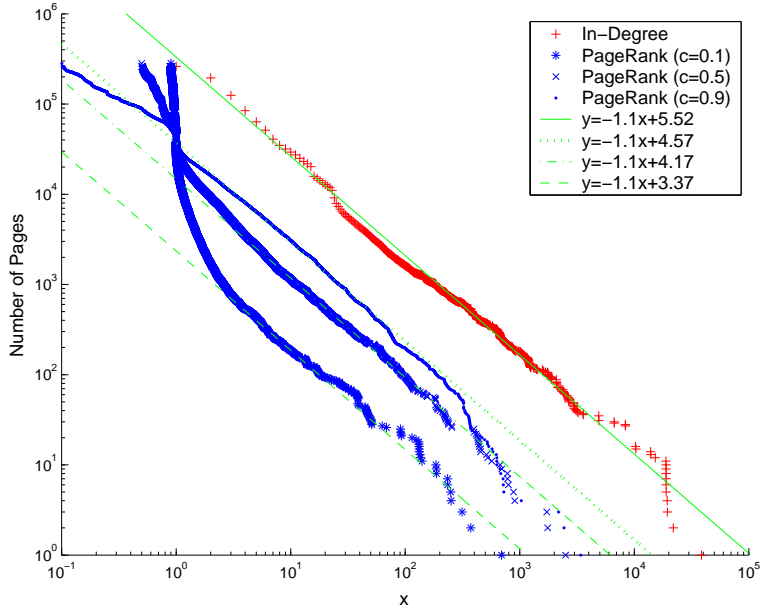}}
        \caption{Plots for the Web data. Number of pages with In-Degree/PageRank greater than $x$  versus $x$ in log-log scale, and the fitted straight lines.}
        \label{web}
\end{figure}

There are several papers, see \cite{Becchetti}, \cite{Fortunato},
\cite{Donato} and \cite{Pandurangan}, that describe similar
experiments for different domains and different number of pages,
and they all confirm that the PageRank and the In-Degree follow
power laws with the same exponent, around 2.1. In Figure~\ref{web}
we show the log-log plots for the In-Degree and the PageRank of
the Stanford Web Data, for different values of the damping factor
($c=0.1,\enskip0.5$ and $0.9$). Clearly, these empirical values of
In-Degree and PageRank constitute parallel straight lines for all
values of the damping factor, provided that the PageRank values
are reasonably large. It was observed in \cite{Becchetti} that in
general, the PageRank depends on the damping factor but the
PageRank of the top 10\% of pages obeys a power law with the same
exponent as the In-Degree, independent on the damping factor. This
is in perfect agreement with our experimental results and the
mathematical model, which is focused on the right tail behavior of
the PageRank distribution.

\begin{figure}[hbt]
            \centering {\epsfxsize=4 in \epsfbox{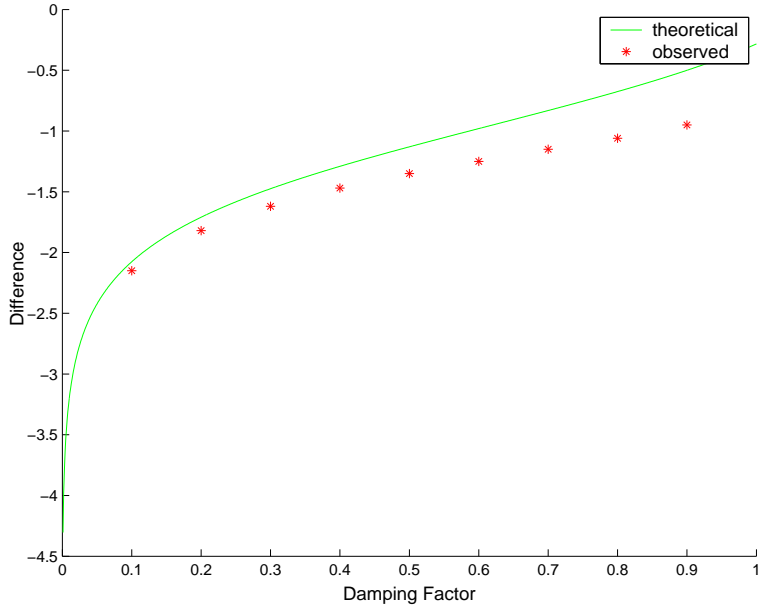}}
        \caption{The theoretical and observed differences between
        logarithmic asymptotics of the In-Degree and the PageRank.}
        \label{diff}
\end{figure}

The calculations based on the maximum likelihood method yield a
slope $-1.1$ for each of the lines, which  verifies that the
In-Degree and PageRank have power laws with the same exponent
$\alpha= 1.1$ (which corresponds to the well known value $2.1$ for
the histogram).  More precisely, we fitted the lines
$y=-1.1x+5.52,\enskip y=-1.1x+4.57,\enskip y=-1.1x+4.17$, and
$y=-1.1x+3.37$ for the plots of the In-Degree and PageRanks with
$c=0,9$, $c=0.5$ and $c=0.1$, respectively.
We also investigated whether Theorem~\ref{t2} correctly predicts
the multiplicative factor
\begin{equation*}
y(c)=\frac{c^{\alpha}}{d^{\alpha}-c^{\alpha}d}.
\end{equation*}
In Figure~\ref{diff} we plotted $\log_{10}(y(c))$ and we compared
it to the observed differences between the logarithms of the
complementary cumulative distribution functions of the PageRank
and the In-Degree, for different values of the damping factor.
Here $d=8.2$ as in the Web data. We see that theoretical and
observed values are remarkably close. Thus, our model not only
allows to prove the similarity in the power law behavior but also
gives a good approximation for the difference between the two
distributions.

The discrepancy between the predicted and observed values of the
multiplicative factor suggests that our model does not capture the
PageRank behavior to the full extent. For instance, the assumption
of the independence of the PageRank of pages that have a common
neighbor may be too strong. We believe however that  the achieved
precision, especially for small values of $c$, is quite good for
our relatively simple stochastic model.

\subsection{Growing Networks} \label{sec:GN}
Growing Networks, introduced by Barab\'asi and Albert
\cite{Barabasi}, now represent a large class of models that are
commonly accepted as a possible scenario of Web growth. In
particular, these models provide a mathematical explanation for
the power law behavior of the In-Degree~\cite{Bollobas}. The
recent studies \cite{Avrachenkov2006}, \cite{Fortunato2} addressed
for the first time the PageRank distribution in Growing Networks.

Growing Network models are characterized by preferential
attachment. This entails that a newly created node connects to the
existing nodes with probabilities that are proportional to the
current In-Degrees of the existing nodes. We simulated a slightly
modified version of this model, where a new link points to a
randomly chosen page with probability $\beta$, and with
probability $1-\beta$ the preferential attachment selection rule
is used. This allows us to tune the exponent of the resulting
power law~\cite{Newman}.

We simulate our Growing Network using Matlab. We start with $d$
nodes and at each step we add a new node that links to $d$ already
existing nodes. To ensure the same number of outgoing links for
all pages, at the end of the simulation, we link the first $d$
nodes to randomly chosen pages. In the example presented below we
set $\beta=0.2$ and obtain a network of $50000$ nodes with
Out-Degree $d=8$.

\begin{figure}[hbt]
            \centering {\epsfxsize=4 in \epsfbox{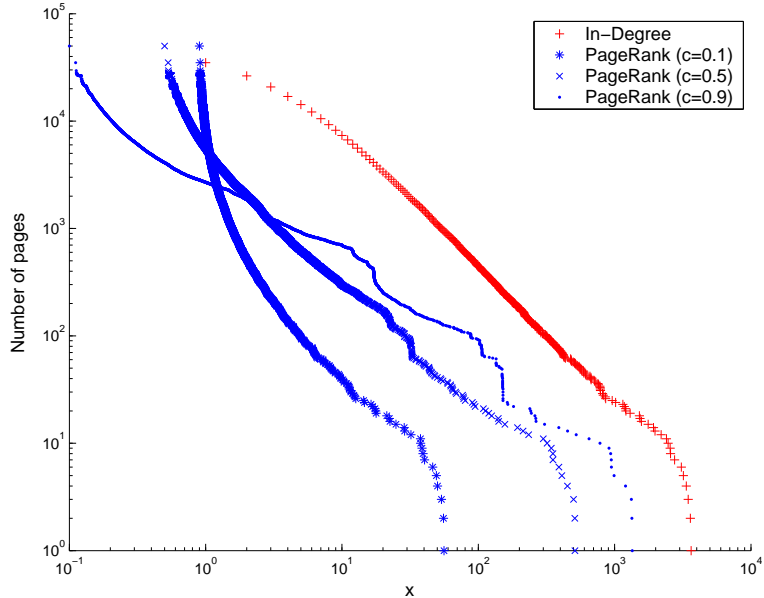}}
        \caption{Plots for the Growing Network model. Number of pages with In-Degree/PageRank greater than $x$  versus $x$ in log-log scale. }
        \label{picturegrnet}
\end{figure}

In Figure~\ref{picturegrnet} we present the numerical data for the
In-Degree and the PageRank in the Growing Network. Clearly, the
Web data from Section~\ref{sec:WD} shows a much better agreement
with our model than the data generated by the preferential
attachment algorithm. In the next section we briefly compare
recent results on the PageRank in Growing Networks to our present
study and we indicate possible directions for further research.

\section{Discussion}

Our model and analysis resulted in the conclusion that the
PageRank and the  In-Degree should follow power laws with the same
exponent. Growing Network models may provide an alternative
explanation~\cite{Avrachenkov2006,Fortunato2}. For instance, in
the recent paper by Avrachenkov and Lebedev~\cite{Avrachenkov2006}
it was shown that the {\em expected} PageRank in Growing Networks
follows a power law with an exponent, which does depend on the
damping factor but equals $\approx 2.08$ for $c=0.85$.  Thus, the
model in \cite{Avrachenkov2006} can also be used to explain the
tail behavior of the PageRank, but it leads to a slightly
different result than our model because in our case the power law
exponent of the PageRank does {\it not} depend on the damping
factor. The reason could be that we focus only on the asymptotics,
whereas \cite{Avrachenkov2006} employs a mean-field approximation.
Indeed, experiments show that the shape of the PageRank
distribution does depend on the damping factor, and thus, it may
affect the average values, whereas the tail behavior remains the
same for all values of $c$.

We emphasise that compared to \cite{Avrachenkov2006,Fortunato2},
our model provides a completely different approach for modelling
the relation between the In-Degree and the PageRank. Specifically,
we do not make any assumption on the underlying Web graph, whereas
\cite{Avrachenkov2006,Fortunato2} choose for the preferential
attachment structure, thus exploiting the fact that this graph
model correctly captures the In-Degree distribution. We believe
that both approaches should be elaborated and used in further
research on the  PageRank distribution.

One of the important innovations in the present work is the
analogy between the PageRank equation and the equation for the
busy period that enables us to apply the techniques from
\cite{MeyerTeugels}. In fact, queueing systems with heavy tails
and in particular the busy period problem allow for a more
sophisticated probabilistic analysis (see e.g. \cite{Zwart}). It
would be interesting to apply these advanced methods to the
problems related to the World Wide Web and PageRank.

Our model definitely lacks the dependencies between the PageRanks
of the pages sharing a common neighbor. Such dependencies must be
present in the Web in particular due to the high clustering of the
Web graph~\cite{Newman} (roughly speaking, clustering means that with high probability, two neighbors of the same page are connected to each other). Thus, in our further research
we could try to include some sort of dependencies in our
stochastic equation. Another natural way to bring our model closer
to the real-life situation is to allow random (heavy-tailed)
Out-Degrees. It would be interesting to investigate in
which ways these new features will affect the PageRank
asymptotics.

\end{document}